\documentclass[a4paper]{amsart}
\usepackage[english]{babel}
\usepackage{amssymb}
\usepackage{amsmath} 
\usepackage{amsthm} 
\usepackage{mathbbol}
\usepackage{graphicx} 
\usepackage{xcolor}
\usepackage{todonotes}
\usepackage{verbatim}

\theoremstyle{plain} 
\newtheorem{thm}{Theorem}  
\newtheorem{prop}[thm]{Proposition}
\newtheorem{lemma}[thm]{Lemma}

\DeclareMathOperator{\Rp}{Re}

\title{On the joint second moment of zeta and its logarithmic derivative}
\author{Alessandro Fazzari}
\address{D\'epartement de math\'ematiques et de statistique, Universit\'e de Montr\'eal. CP 6128, succ. Centre-ville. Montreal, QC H3C 3J7, Canada}
\email{alessandro.fazzari@umontreal.ca}
 
\begin{document}
\maketitle

%\begin{abstract}
%Assuming the Riemann hypothesis, we study the second moment of zeta times its logarithmic derivative, shifted away from the half line by $a/\log T$. 
%This can be seen as weighted version of a result by Goldston, Gonek and Montgomery, where the authors consider the second moment of the log-derivative of $\zeta$. Here,  we provide an upper and a lower bound for the second moment of $\frac{\zeta'}{\zeta}(\frac{1}{2}+\frac{a}{\log T}+it)$, where the average is tilted by $|\zeta(\frac{1}{2}+it)|^2$.
%As a corollary, we provide an upper bound for the variance on short intervals of the sum of the arithmetic function given by the von Mangoldt convolved with a shift.
%\end{abstract}

\begin{abstract}
Assuming the Riemann Hypothesis, Goldston, Gonek and Montgomery \cite{GGM} studied the second moment of the log-derivative of $\zeta$ shifted away from the half line by $a/\log T$, and its connection with the pair correlation conjecture. 
In this paper we consider a weighted version of this problem where the average is tilted by $|\zeta(\frac{1}{2}+it)|^2$.
We provide an upper and a lower bound for the second moment of zeta times its logarithmic derivative, $a/\log T$ away from the critical line. 
%As a corollary, we deduce an upper bound for the variance on short intervals of the sum of the arithmetic function given by the von Mangoldt convolved with a shift.
\end{abstract}

\section{Introduction}\label{SectionMainResult}

In 1973, Montgomery \cite{Montgomery} studied the distribution of pairs of zeros of the Riemann zeta function. Assuming the Riemann Hypothesis (RH) and denoting $\rho=\frac{1}{2}+i\gamma$ with $\gamma\in\mathbb R$ the nontrivial zeros of $\zeta$, Montgomery conjectured that for any $a<b$
\begin{equation}\label{PCC}
\frac{1}{N(T)}\sum_{\substack{0<\gamma,\gamma'\leq T \\ \frac{2\pi a}{\log T}<\gamma-\gamma'<\frac{2\pi b}{\log T}}}1
\sim \int_a^b \left(1-\frac{\sin(\pi t)^2}{(\pi t)^2}\right)dt +\delta(a,b)
\end{equation}
as $T\to \infty$, where $N(T)=\symbol{35}\{0<\gamma\leq T\}\sim \frac{T\log T}{2\pi}$ and $\delta(a,b)=1$ if $0\in[a,b]$ and 0 otherwise. Equation \eqref{PCC}, known as the pair correlation conjecture, suggests that the presence of a given zero $\gamma$ at height $T$ repells zeros $\gamma'$ nearby. 
%Namely, the probability of finding a zero $\gamma'$ such that $0<|\gamma-\gamma'|\leq 2\pi x/\log T$ tends to $\int_{-x}^{x}  (1-\frac{\sin(\pi t)^2}{(\pi t)^2})dt$, as $T\to\infty $. %$\asymp x^3$, as $x\to0$.
Namely, the probability of finding a zero $\gamma'$ that is less than $2\pi x/\log T$ away from another zero $\gamma$ is given by $\int_{-x}^{x}  (1-\frac{\sin(\pi t)^2}{(\pi t)^2})dt$. %$\asymp x^3$, as $x\to0$.

%$$\mathbb P \bigg[\gamma': 0<|\gamma-\gamma'|\leq \frac{2\pi x}{\log T}\bigg] \sim \int_{-x}^{x}  \left(1-\frac{\sin(\pi t)^2}{(\pi t)^2}\right)dt.$$
%If $x\to 0$, the above goes to zero at the same rate as $x^3$, predicting repulsion between zeros.

In 2001, Goldston, Gonek and Montgomery \cite{GGM} proved (under RH) that the pair correlation conjecture is equivalent to the asymptotic formula
\begin{equation}\label{8dec.1} 
\int_{T}^{2T} \left|\frac{\zeta'}{\zeta}\left(\frac{1}{2}+\frac{a}{\log T}+it\right) \right|^2 dt \sim \left(\frac{1-e^{-2a}}{4a^2}\right)T(\log T)^2 
\end{equation}
as $T\to\infty$ and for any $a>0$. %, such that $T^{-1}(\log T)^3\ll a \ll 1$. 
Note that, if $a=o(1)$ as $T\to\infty$, the function $(1-e^{-2a})/4a^2$ is asymptotic to $\frac{1}{2a}$ and therefore the above quantity explodes as $a$ goes to $0$. This is due to the effect of the nontrivial zeros of zeta that make the integral diverge as we approach the critical line.

Furthermore, under RH only, they studied the left-hand side of \eqref{8dec.1}, providing upper and lower bound when $(\log T)^3/T \leq a = o(1)$ as $T\to\infty$; they showed that \cite[Corollary 3]{GGM}
\begin{equation}\label{Cor3GGM} 
\frac{1-o(1)}{2a}
\leq  \frac{1}{T(\log T)^2}\int_{T}^{2T} \left|\frac{\zeta'}{\zeta}\left(\frac{1}{2}+\frac{a}{\log T}+it\right) \right|^2 dt 
\leq \frac{\frac{4}{3}+o(1)}{2a}.
\end{equation}
In this context, we also mention a recent paper by Carneiro, Chirre, Chandee and Milinovich \cite{CCCM}. They improved the bounds in \eqref{Cor3GGM} by introducing new connections to certain extremal problems in Fourier analysis.

This paper concerns a weighted version of the problem described above. The main motivation for this approach is given by its connection to the interplay between zeros and large values of zeta.  Following a similar philosophy as in \cite{3., 4.}, we are interested in how the presence of a large value of zeta influences the distribution of zeros closeby. To be concrete, suppose that we want to look at how pairs of zeros are distributed around those points $t$ where zeta is large, say of size about $(\log T)^k$. Note that, on RH, we know that these large values are responsible for the $2k$-th moment of zeta (see \cite{Sound}). Therefore, it seems natural to weight the measure of the integral in \eqref{8dec.1} by the factor $|\zeta(\frac{1}{2}+it)|^{2k}$. These weights already appeared in \cite{1., 2., BELP}, concerning weighted central limit theorems for central values of $L$-functions. 

For example, in the case $k=1$, the problem is that of studying the quantity
$$ I_{|\zeta|^2}( a ;T) := \int_{T}^{2T} \left|\frac{\zeta'}{\zeta}\left(\frac{1}{2}+\frac{ a }{\log T}+it\right) \right|^2 |\zeta(\tfrac{1}{2}+it)|^2 dt$$
for $a>0$ and as $T\to\infty$. The quantity $I_{|\zeta|^2}(a;T)$ is sensitive to the zeros and the large values simultaneously. Indeed, in the integral above, the first factor encodes information about the distribution of pair of zeros, while the second factor tends to concentrate the integral around those values of $t$ such that $|\zeta(\frac{1}{2}+it)|\asymp(\log T)^{1+o(1)}$.

Using the heuristic machinery of \cite{CFKRS}, we can predict the expected asymptotic formula for $I_{|\zeta|^2}( a ;T)$. More precisely, assuming the ratio conjecture \cite[Conjecture 5.1]{CFZ} with four copies of zeta at the numerator and two at the denominator, one can show that for $a\ll 1$
%\begin{thm}\label{conditionalthm}
%Assume RH and Conjecture 5.1 of \cite{CFZ} in the case $K=L=2$ and $Q=R=1$. Then, for $a\ll 1$, we have
\begin{equation}\label{expectation} \notag
\frac{I_{|\zeta|^2}( a ;T) }{T(\log T)^3} = f(a)+ O\bigg(\frac{1}{\log T}\bigg)
\end{equation}
with
$$f(a):=\frac{(5a - 8) + 12 e^{-a} - ( a + 4)e^{-2a} }{4a^3}.$$
%\end{thm}

We notice that, as $a\to 0$, the main term above is
\begin{equation}\label{asympf} 
f(a) = \frac{1}{3}-\frac{5}{24}a + \frac{3}{40}a^2 - \frac{13}{720}a^3 %+ \frac{13}{5040}a^4+ \frac{1}{13440}a^5 - \frac{1}{630}a^6 +\dots
+O(a^4). \end{equation}
Hence, $I_{|\zeta|^2}( a ,T)$ converges to $\frac{1}{3}T(\log T)^3$, i.e. the main term of the second moment of $\zeta'(\frac{1}{2}+it)$, unlike the classical quantity in \eqref{8dec.1}, which diverges for $a\to0$. This is due to the weight $|\zeta(\tfrac{1}{2}+it)|^2$ in the integral over $t$, which kills the contribution from zeros that are responsible for the divergence of the moments of the logarithmic derivative of zeta.
%Theorem \ref{conditionalthm} immediately follows from the assumption of the ratio conjecture, since
%$$ I_{|\zeta|^2}(a,T) = \partial_E\partial_F \bigg[\int_{T}^{2T} \frac{\zeta(\tfrac{1}{2}+A+it)\zeta(\tfrac{1}{2}+B-it)\zeta(\tfrac{1}{2}+E+it)\zeta(\tfrac{1}{2}+F-it)}{\zeta(\tfrac{1}{2}+C+it)\zeta(\tfrac{1}{2}+D-it)} dt \bigg]_{\substack{E=C= a /\log T \\ F=D= a /\log T \\ A=B=0 }}.$$ 

Without any assumptions on the moments of zeta, the main result of this paper provides upper and lower bound for $I_{|\zeta|^2}(a,T)$. 

%\textcolor{red}{
%\begin{thm}
%Assume RH. If $(\log T)^{-1/4+\varepsilon}\ll a \ll 1$ for some $\varepsilon>0$, we have 
%$$  \frac{ 4 \sinh(a) + 2 (a - 5) \cosh(a) -a^3 + 5 a^2 - 6 a +10}{2 a^3}+ O((\log T)^{-\varepsilon}) \leq \frac{I_{|\zeta|^2}( a ;T) }{T(\log T)^3} \leq  \frac{e^{-a}(-4a^2-4a-2)+2e^a}{(2a)^3} + O((\log T)^{-1}).$$
%In particular, if $a=o(1)$, the above gives
%$$ \bigg(\frac{1}{3} - \frac{5 a}{24}\bigg)+ O((\log T)^{-\varepsilon}) \leq \frac{I_{|\zeta|^2}( a ;T) }{T(\log T)^3} \leq  \bigg(\frac{1}{3} - \frac{a}{6} + \frac{a^2}{15}\bigg) + O((\log T)^{-1}).$$
%\end{thm}
%}
%
%\textcolor{red}{
%\begin{thm}
%Assume RH. If $(\log T)^{-1/4+\varepsilon}\ll a \ll 1$ for some $\varepsilon>0$, we have 
%$$  \frac{ 4 \sinh(a) + 2 (a - 5) \cosh(a) -a^3 + 5 a^2 - 6 a +10}{2 a^3} \leq \frac{I_{|\zeta|^2}( a ;T) }{T(\log T)^3} + O\bigg(\frac{1}{(\log T)^{\varepsilon}}\bigg) \leq  \frac{e^{-a}(-4a^2-4a-2)+2e^a}{(2a)^3}.$$
%In particular, if $a=o(1)$, the above gives
%$$ \bigg(\frac{1}{3} - \frac{5 a}{24}\bigg) \leq \frac{I_{|\zeta|^2}( a ;T) }{T(\log T)^3}  + O\bigg(\frac{1}{(\log T)^{\varepsilon}}\bigg) \leq  \bigg(\frac{1}{3} - \frac{a}{6} + \frac{a^2}{15}\bigg) .$$
%\end{thm}
%}

%\textcolor{red}{
\begin{thm}\label{mainthm}
Assume RH. Let
\begin{equation}\begin{split}\notag 
f_1(a):&=\frac{ - a^3 +5 a^2 -10 a - 10 (e^{-a}-1)}{2 a^3} 
= \frac{1}{3} - \frac{5 a}{24} + \frac{a^2}{24} - \frac{a^3}{144} + O(a^4), \quad \text{as } a\to0 \\
f_2(a):&= \frac{e^a-e^{-a}(2a^2+2a+1)}{4a^3} 
=\frac{1}{3} - \frac{a}{6} + \frac{a^2}{15} - \frac{a^3}{60} + O(a^4), \quad \text{as } a\to0.
\end{split}\end{equation}
If $(\log T)^{-1/4+\varepsilon}\ll a \ll 1$ for some $0<\varepsilon\leq\frac{1}{4}$, then we have 
%$$   f_1(a)T(\log T)^3+ O(T(\log T)^{3-\varepsilon}) \leq I_{|\zeta|^2}( a ;T) \leq f_2(a)T(\log T)^3 + O(T(\log T)^{2}).$$
%$$   f_1(a)+ O((\log T)^{-\varepsilon}) \leq \frac{ I_{|\zeta|^2}( a ;T)}{T(\log T)^3} \leq f_2(a) + O((\log T)^{-1}).$$
$$   f_1(a)+ O\bigg(\frac{1}{(\log T)^\varepsilon}\bigg) \leq \frac{ I_{|\zeta|^2}( a ;T)}{T(\log T)^3} \leq f_2(a) + O\bigg(\frac{1}{\log T}\bigg).$$
In particular, if $a=o(1)$, the above gives
$$ \bigg(\frac{1}{3} - \frac{5 a}{24}\bigg)+ O\bigg(\frac{1}{(\log T)^\varepsilon}\bigg) \leq \frac{I_{|\zeta|^2}( a ;T) }{T(\log T)^3} \leq  \bigg(\frac{1}{3} - \frac{a}{6} + \frac{a^2}{15}\bigg) + O\bigg(\frac{1}{\log T}\bigg).$$
\end{thm}
%}

We notice that the bounds provided by the theorem reflect the expected order of magnitude, of size $T(\log T)^3$.  In terms of the constant, neither the lower bound nor the upper bound are sharp, since $f_1(a)<f(a)<f_2(a)$ for $a>0$. However, when $a\to0$, both the upper and lower bounds tend to $\frac{1}{3}$, in accordance with \eqref{asympf}. Moreover, the lower bound captures also the second term $-\frac{5}{24}a$ term appearing in \eqref{asympf}. 

This result should be compared to \eqref{Cor3GGM}. Indeed, Theorem \ref{mainthm} can be seen as a generalization of the result by Goldston, Gonek and Montgomery, when the integral over $t$ is tilted by $|\zeta(\frac{1}{2}+it)|^2$. In both cases, the lower bound is more accurate than the upper bound, as it matches the expected asymptotic behavior as $a\to 0$.

The paper is structured as follows. In Section \ref{SectionLowerBound}, we prove the lower bound appearing in Theorem \ref{mainthm}. First we apply the explicit formula to express $|\frac{\zeta'}{\zeta}(\frac{1}{2}+\frac{a}{\log T}+it)|^2$ in terms of a sum over pairs of zeros $\sum_{\gamma,\gamma'}g(\gamma,\gamma')$ for an explicit positive function $g$. Then, we lower bound the double sum by the diagonal contribution $\gamma=\gamma'$. This diagonal contribution can be evaluated asymptotically on weighted average with respect to the measure $|\zeta(\frac{1}{2}+it)|^2dt$ due to Gonek's formula for the discrete second moment of zeta on RH. For the upper bound in Theorem \ref{mainthm}, provided in Section \ref{SectionUpperBound}, we use a result by Soundararajan that (on RH) bounds the ratio $\zeta(\frac{1}{2}+it)/\zeta(\frac{1}{2}+\frac{a}{\log T}+it)$, and the well-known formula for the second moment of the derivative of zeta. %Finally, Theorem \ref{thmprimes} is proven in Section \ref{Sectionprimes} by using the upper bound for $I_{|\zeta|^2}(a,T)$.

\medskip
\textbf{Acknowledgments}. The author would like to thank Steve Lester for suggesting the problem and for several inspiring conversations about the topic. Also, thanks to Sandro Bettin and Brian Conrey for many helpful discussions. This work was supported by the Fonds de recherche du Qu\'ebec - Nature et technologies, Projet de recherche en \'equipe 300951 and the FRG grant DMS 1854398. The author is a member of the INdAM group GNAMPA.

\section{Lower bound}\label{SectionLowerBound}

In this section we prove the lower bound for $I_{|\zeta|^2}(a;T)$.
%\begin{prop}\label{lowerbound}
%Assume RH. For $a\gg\frac{\log T}{T}$, as $T\to\infty$, we have
%\begin{equation}\begin{split}\notag I_{|\zeta|^2}( a ;T) \geq T(\log T)^3 &\frac{ - a^3 +5 a^2 -10 a - 10 (e^{-a}-1)}{2 a^3}  
%+ O(a^2T) + O\bigg(\frac{T(\log T)^{11/4}}{a}\bigg) + O(T(\log T)^2). 
%\end{split}\end{equation}
%In particular, if $(\log T)^{-1/4+\varepsilon}\ll a \ll 1$ for some $\varepsilon>0$, then
%\begin{equation}\notag I_{|\zeta|^2}( a ;T) \geq T(\log T)^3 \frac{ - a^3 +5 a^2 -10 a - 10 (e^{-a}-1)}{2 a^3}   + O(T(\log T)^{3-\varepsilon}). \end{equation}
%If also $a=o(1)$, then
%\begin{equation}\notag I_{|\zeta|^2}( a ;T) \geq T(\log T)^3 \bigg(\frac{1}{3} - \frac{5 a}{24}\bigg)+ O(T(\log T)^{3-\varepsilon}). 
%\end{equation}
%\end{prop}

\begin{prop}\label{lowerbound}
Assume RH. For $0< a \ll 1$, as $T\to\infty$, we have
\begin{equation}\begin{split}\notag I_{|\zeta|^2}( a ;T) \geq T(\log T)^3 &\frac{ - a^3 +5 a^2 -10 a - 10 (e^{-a}-1)}{2 a^3}  
+ O\bigg(\frac{T(\log T)^{11/4}}{a}\bigg) + O\bigg(\frac{T(\log T)^2}{a^2}\bigg) .%+ O(T(\log T)^2). 
\end{split}\end{equation}
In particular, if $(\log T)^{-1/4+\varepsilon}\ll a \ll 1$ for some $0<\varepsilon\leq\frac{1}{4}$, then
%\begin{equation}\notag I_{|\zeta|^2}( a ;T) \geq T(\log T)^3 \frac{ - a^3 +5 a^2 -10 a - 10 (e^{-a}-1)}{2 a^3}   + O(T(\log T)^{3-\varepsilon}). \end{equation}
the error term above is $\ll T(\log T)^{3-\varepsilon}$.\quad\quad\quad
If also $a=o(1)$, then
\begin{equation}\notag I_{|\zeta|^2}( a ;T) \geq T(\log T)^3 \bigg(\frac{1}{3} - \frac{5 a}{24}\bigg)+ O(T(\log T)^{3-\varepsilon}). 
\end{equation}
\end{prop}
%
%We begin with an easy lemma that we will use several times in the following.
%\begin{lemma}\label{lemmaByParts}
%For any differentiable function $f$, we have
%\begin{equation}\begin{split}\notag 
%\int_T^{2T} |\zeta(\tfrac{1}{2}+it)|^2 f(t) dt 
%= \int_T^{2T} (\log t + 2\gamma-\log(2\pi)) f(t) dt &+ O(\sqrt{T}|f(2T)|) 
%\\&+O(\sqrt{T}|f(T)|)  + O\bigg(\int_T^{2T} \sqrt{t}|f'(t)|dt\bigg). 
%\end{split}\end{equation}
%\end{lemma}
%\proof It follows immediately by integration by parts and the asymptotic formula for the second moment of zeta $\int_0^t |\zeta(\frac{1}{2}+iy)|^2dy = t\log t + (2\gamma - 1 -\log(2\pi)) t +O(\sqrt{t})$. \endproof
%
The first step towards Proposition \ref{lowerbound} consists of expressing the logarithmic derivative of zeta in $I_{|\zeta|^2}(a;T)$ in terms of a sum over zeros. This can be done by an application of the explicit formula, which gives the following.

\begin{lemma}\label{lemmaSpezzamentoI}
Assume RH. For $0<a\ll 1$, we have
$$ I_{|\zeta|^2}( a ;T) 
=  T(\log T)^3 \frac{(a-2)(2e^{-a}-2-a^2+2a)}{2 a^3} + 2 \mathcal S_{|\zeta|^2}\left(\frac{1}{2}+\frac{a}{\log T},T\right) + O\bigg(\frac{T(\log T)^2}{a^2}\bigg), $$
as $T\to\infty$, where
$$ \mathcal S_{|\zeta|^2}(\sigma;T):= \int_T^{2T} \left(\sum_\gamma\frac{\sigma-\frac{1}{2}}{(\sigma-\frac{1}{2})^2+(t-\gamma)^2} \right)^2|\zeta(\tfrac{1}{2}+it)|^2dt. $$
\end{lemma}

\proof
We follow closely the strategy of the proof of Theorem 1 in \cite{GGM}. Being $|w|^2=2(\Rp w)^2-\Rp(w^2)$, we have (denoting $\sigma=\frac{1}{2}+\frac{a}{\log T}$)
\begin{equation}\label{14oct.1}
I_{|\zeta|^2}( a ;T) 
=  2\int_{T}^{2T} \left(\Rp\tfrac{\zeta'}{\zeta}\left(\sigma+it\right)\right)^2 |\zeta(\tfrac{1}{2}+it)|^2 dt
-\Rp \int_{T}^{2T} \left(\tfrac{\zeta'}{\zeta}\left(\sigma+it\right)\right)^2 |\zeta(\tfrac{1}{2}+it)|^2 dt.
\end{equation}
To analyze the first term, we recall that, denoting 
%\begin{equation}\begin{split}\notag
%\mathcal L(\sigma;T) :&= \int_T^{2T} \left(\Rp\tfrac{\zeta'}{\zeta}(\sigma+it)+\tfrac{1}{2}\log\tfrac{t}{2\pi}+O(\tfrac{1}{t})\right)^2 |\zeta(\tfrac{1}{2}+it)|^2dt \\
%\mathcal S_{|\zeta|^2}(\sigma;T):&= \int_T^{2T} \left(\sum_\gamma\frac{\sigma-\frac{1}{2}}{(\sigma-\frac{1}{2})^2+(t-\gamma)^2} \right)^2|\zeta(\tfrac{1}{2}+it)|^2dt
%\end{split}\end{equation}
\begin{equation}\notag
\mathcal L(\sigma;T) := \int_T^{2T} \left(\Rp\tfrac{\zeta'}{\zeta}(\sigma+it)+\tfrac{1}{2}\log\tfrac{t}{2\pi}+O(\tfrac{1}{t})\right)^2 |\zeta(\tfrac{1}{2}+it)|^2dt
\end{equation}
we have (see \cite{GGM} p.112, above Equation (2.4))
\begin{equation}\label{L=R}
\mathcal L(\sigma;T) = \mathcal S_{|\zeta|^2}(\sigma;T).
\end{equation}
Expanding the square and bounding (by the well-known bound $ \frac{\zeta'}{\zeta}(u+iT) \ll \frac{\log T}{u-\frac{1}{2}}$ for $\frac{1}{2}<u\leq 2$) the error terms, one gets
\begin{equation}\begin{split}\label{22may.01}
\mathcal L(\sigma;T) &= \int_T^{2T} \left(\Rp\tfrac{\zeta'}{\zeta}(\sigma+it)\right)^2 |\zeta(\tfrac{1}{2}+it)|^2dt + \int_T^{2T} \tfrac{1}{4}(\log\tfrac{t}{2\pi})^2 |\zeta(\tfrac{1}{2}+it)|^2dt \\
&\hspace{2.5cm}+ \Rp\int_T^{2T} \tfrac{\zeta'}{\zeta}(\sigma+it)\log\tfrac{t}{2\pi} |\zeta(\tfrac{1}{2}+it)|^2dt + O\left( \frac{(\log T)^3}{a} \right).
\end{split}\end{equation}
Moreover, being $\log\frac{t}{2\pi}=\log T+O(1)$ for $T\leq t\leq 2T$, the asymptotic formula for the second moment of zeta yields
\begin{equation}\label{22may.02} 
 \int_T^{2T} \tfrac{1}{4}(\log\tfrac{t}{2\pi})^2 |\zeta(\tfrac{1}{2}+it)|^2dt = \tfrac{1}{4}T(\log T)^3 +O(T(\log T)^2).
\end{equation}
Plugging \eqref{22may.01} and \eqref{22may.02} into \eqref{L=R}, we obtain
\begin{equation}\begin{split}\notag
 \int_T^{2T}& \left(\Rp\tfrac{\zeta'}{\zeta}(\sigma+it)\right)^2 |\zeta(\tfrac{1}{2}+it)|^2dt 
 =\mathcal S_{|\zeta|^2}(\sigma,T) - \tfrac{1}{4}T(\log T)^3  \\
&\quad
- \Rp\int_T^{2T} \tfrac{\zeta'}{\zeta}(\sigma+it)\log\tfrac{t}{2\pi} |\zeta(\tfrac{1}{2}+it)|^2dt + O(T(\log T)^2) + O\left( \frac{(\log T)^3}{a} \right).
\end{split}\end{equation}
Then, Equation \eqref{14oct.1} reads
\begin{equation}\begin{split}\label{splitting}\notag
I_{|\zeta|^2}( a ;T) 
&= -2\Rp I_{|\zeta|^2}^{(1)}(a;T)  - \tfrac{1}{2}T(\log T)^3 -\Rp I_{|\zeta|^2}^{(2)}(a;T)\\
&\quad\quad+2 \mathcal S_{|\zeta|^2}(\sigma;T) + O(T(\log T)^2) + O\left( \frac{(\log T)^3}{a} \right)
\end{split}\end{equation}
with
\begin{equation}\begin{split}\label{I1}
I_{|\zeta|^2}^{(1)}( a ;T) :&=  \int_T^{2T} \tfrac{\zeta'}{\zeta}(\sigma+it)\log\tfrac{t}{2\pi} |\zeta(\tfrac{1}{2}+it)|^2dt
= T(\log T)^3 \frac{1-a-e^{-a}}{a^2} + O\bigg(\frac{T(\log T)^2}{a}\bigg),
\end{split}\end{equation}
\begin{equation}\begin{split}\label{I2}
I_{|\zeta|^2}^{(2)}( a ;T) :&= \int_{T}^{2T} \left(\tfrac{\zeta'}{\zeta}\left(\sigma+it\right)\right)^2 |\zeta(\tfrac{1}{2}+it)|^2 dt
= T(\log T)^3  \frac{ a-2 + (a+2)e^{-a} }{a^3} + O\bigg(\frac{T(\log T)^2}{a^2}\bigg),
\end{split}\end{equation}
that we will prove below. This gives the claim.
%by \cite[Lemma 2]{ConreySimpleZeros}).
%\todo[inline]{add proof of these. Anyway I think I know how to do it. Just approximate the logarithmic derivative of zeta with a D-poly (see e.g. (3.2) in \lq\lq The distribution of the logarithmic derivative of the Riemann zeta-function\rq\rq, by Steve) and then use essentially the same computation as in my paper with Sandro, i.e. twisted second moment $\int P\times|\zeta|^2$.}
\endproof

\proof (of Equation \eqref{I1})
For starters, we note that it suffices to show that 
\begin{equation}\label{23may.01}
\int_1^t \frac{\zeta'}{\zeta}(\sigma+iy)|\zeta(\tfrac{1}{2}+iy)|^2 dy = t(\log t)^2 \frac{1-a-e^{-a}}{a^2} + O\bigg(\frac{t\log t}{a}\bigg)
\end{equation}
for $T\leq t\leq 2T$, $\sigma=\frac{1}{2}+\frac{a}{\log T}$. Indeed, Equation \eqref{I1} follows by integration by parts; denoting $h(a)=\frac{1-a-e^{-a}}{a^2}$, Equation \eqref{23may.01} yields
\begin{equation}\begin{split}\notag
I_{|\zeta|^2}(a;T) 
&= \bigg[  \log\tfrac{t}{2\pi}  \int_1^t \frac{\zeta'}{\zeta}(\sigma+iy)|\zeta(\tfrac{1}{2}+iy)|^2 dy \bigg]_T^{2T} 
- \int_T^{2T} \int_1^t \frac{\zeta'}{\zeta}(\sigma+iy)|\zeta(\tfrac{1}{2}+iy)|^2 dy \frac{d}{dt}\bigg[\log\frac{t}{2\pi}\bigg]dt\\
&= h(a) \bigg[  t(\log t)^2 \log\tfrac{t}{2\pi} \bigg]_T^{2T} 
- h(a) \int_T^{2T}t(\log t)^2  \frac{d}{dt}\bigg[\log\frac{t}{2\pi}\bigg]dt
+ O\bigg(\frac{T(\log T)^2}{a}\bigg)\\
&= %g(a) \bigg[  t(\log t)^2 \log\tfrac{t}{2\pi} \bigg]_T^{2T} 
%- g(a) \bigg[t(\log t)^2 \log \tfrac{t}{2\pi}\bigg]_T^{2T}+
 h(a) \int_T^{2T}((\log t)^2+2\log t) \log\tfrac{t}{2\pi} dt
+ O\bigg(\frac{T(\log T)^2}{a}\bigg)
=  h(a)T(\log T)^3
+ O\bigg(\frac{T(\log T)^2}{a}\bigg).
\end{split}\end{equation}
We now prove \eqref{23may.01}. Using the functional equation of the Riemann zeta function 
$$\zeta(1-s) = \chi(1-s)\zeta(s), \quad \text{ with } 
%\chi(s) = \pi^{s-\frac{1}{2}}\frac{\Gamma(\frac{1-s}{2})}{\Gamma(\frac{s}{2})},
% or maybe we want 
\chi(1-s)=\chi(s)^{-1}=2(2\pi)^{-s}\Gamma(s)\cos(\tfrac{\pi s}{2}),
$$
we write 
\begin{equation}\notag
\int_1^t \frac{\zeta'}{\zeta}(\sigma+iy)|\zeta(\tfrac{1}{2}+iy)|^2 dy
=  \int_{\frac{1}{2}+i}^{\frac{1}{2}+it} \frac{\zeta'}{\zeta}\bigg(s+\frac{a}{\log T}\bigg)\zeta(s)^2\chi(1-s)\frac{ds}{i} .
\end{equation}
We shift the integral to the right, on the line $c=1+\frac{1}{\log t }$. The horizontal integrals can be bounded by the Lindel\"of Hypothesis (note that we are assuming RH) and the standard estimates $\chi(\sigma+it)\ll |t|^{\frac{1}{2}-\sigma}$ and  $\frac{\zeta'}{\zeta}(\sigma+it) \ll \frac{\log t}{\sigma-\frac{1}{2}}$. Doing so, we get 
\begin{equation}\begin{split}\notag
\int_1^t \frac{\zeta'}{\zeta}(\sigma+iy)|\zeta(\tfrac{1}{2}+iy)|^2 dy
%&=  2\pi\log T \frac{1}{2\pi i} \int_{c+iT}^{c+2iT} \chi(1-s) \frac{\zeta'}{\zeta}\bigg(s+\frac{a}{\log T}\bigg)\zeta(s)^2ds + O(T(\log T)^2) \\
&=  \frac{2\pi}{2\pi i} \int_{c+i}^{c+it} A(s) \chi(1-s) ds + O\bigg(\frac{t^{\frac{1}{2}+\varepsilon}}{a}\bigg)
\end{split}\end{equation}
for any $\varepsilon>0$, with
$$ A(s) =  \frac{\zeta'}{\zeta}\bigg(s+\frac{a}{\log T}\bigg)\zeta(s)^2
= \sum_{n=1}^{\infty}\frac{g(n)}{n^s}, \quad g(n) = \Big[\Big(-\Lambda(\cdot) \; \cdot^{-a/\log T}\Big)*\tau\Big](n) $$
where $\tau$ denotes the divisor function.
Now we apply Lemma 2 in \cite{ConreySimpleZeros} (with $B(s)=1$), obtaining
\begin{equation}\begin{split}\label{13febr.01}
\int_1^t \frac{\zeta'}{\zeta}(\sigma+iy)|\zeta(\tfrac{1}{2}+iy)|^2 dy
&=  2\pi \sum_{n\leq t/2\pi}g(n) + O\bigg(\frac{t^{\frac{1}{2}+\varepsilon}}{a}\bigg) .
\end{split}\end{equation}
Moreover, by Dirichlet's theorem
\begin{equation}\begin{split}\notag
\sum_{n\leq t/2\pi}g(n) 
%&= - \sum_{n\leq T/2\pi} \sum_{d|n} \Lambda(d)d^{-a/\log T}\tau\bigg(\frac{n}{d}\bigg)
%= - \sum_{d\leq T/2\pi} \Lambda(d)d^{-a/\log T} \sum_{\substack{n\leq T/2\pi \\ n\equiv 0\; (d)}}\tau\bigg(\frac{n}{d}\bigg) \\
& = - \sum_{d\leq t/2\pi} \Lambda(d)d^{-a/\log T} \sum_{n\leq t/(2\pi d)}\tau(n)
= - \sum_{d\leq t/2\pi} \Lambda(d)d^{-a/\log T} \bigg(\frac{t}{2\pi d}\log\frac{t}{2\pi d} + O\bigg(\frac{t}{d}\bigg)\bigg)\\
&= - \frac{t\log t}{2\pi} \bigg(1+O\bigg(\frac{1}{\log t}\bigg)\bigg) \sum_{d\leq t/2\pi} \frac{\Lambda(d)d^{-a/\log T}}{d} 
+ \frac{t}{2\pi} \sum_{d\leq t/2\pi} \frac{\Lambda(d)d^{-a/\log T}\log d}{d} .
\end{split}\end{equation}
Using the stardard asymptotic formulae $\sum_{d\leq x}\Lambda(d)/d = \log t+O(1)$ and $\sum_{d\leq x}\Lambda(d)\log d/d = (\log t)^2/2+O(\log t)$, by partial summation one can easily see that
\begin{equation}\begin{split}\notag
\sum_{d\leq t/2\pi} \frac{\Lambda(d)d^{-a/\log T}}{d} &= \int_1^{t/2\pi} y^{-1-a/\log T}dy +O(1) 
=  \log t\frac{1-e^{-a}}{a} + O(1)\\
\sum_{d\leq t/2\pi} \frac{\Lambda(d)d^{-a/\log T}\log d}{d} & =\int_1^{t/2\pi}  y^{-a/\log T}\frac{\log y}{y}dy  +O(\log t) = (\log t)^2 \frac{1-e^{-a}-ae^{-a}}{a^2} + O\bigg(\frac{\log t}{a}\bigg).
\end{split}\end{equation}
Above, we used the fact that $\log \frac{t}{2\pi} +O(1) = \log t = \log T +O(1)$ for $T\leq t \leq 2T$ and then $(\frac{t}{2\pi})^{-a/\log T} = e^{-a}+O(a/\log T)$. Therefore,
\begin{equation}\begin{split}\notag
\sum_{n\leq t/2\pi}g(n) 
&= -\frac{t(\log t)^2}{2\pi}\frac{1-e^{-a}}{a} + \frac{t(\log t)^2}{2\pi}  \frac{1-e^{-a}-ae^{-a}}{a^2} + O\bigg(\frac{t\log t}{a}\bigg)\\
%&= -\frac{T(\log T)^2}{2\pi}\frac{1-e^{-a}}{a} + \frac{T(\log T)^2}{2\pi}\frac{-e^{-a}-ae^{-a}+1}{a^2} + O(T\log T)\\
&= \frac{t(\log t)^2}{2\pi} \frac{1 - a - e^{-a}}{a^2} + O\bigg(\frac{t\log t}{a}\bigg).
\end{split}\end{equation}
Plugging this into \eqref{13febr.01} we get \eqref{23may.01}.
\endproof

\proof (of Equation \eqref{I2})
Very similarly to what we did above, we write
\begin{equation}\begin{split}\notag
I_{|\zeta|^2}^{(2)}( a ;T) 
%&=  \log T \int_T^{2T} \bigg(\frac{\zeta'}{\zeta}\bigg(\frac{1}{2}+\frac{a}{\log T}+it\bigg)\bigg)^2|\zeta(\tfrac{1}{2}+it)|^2dt \\
&=  \int_{\frac{1}{2}+iT}^{\frac{1}{2}+2iT} \bigg(\frac{\zeta'}{\zeta}\bigg(s+\frac{a}{\log T}\bigg)\bigg)^2\zeta(s)^2\chi(1-s)\frac{ds}{i}\\
&=\frac{2\pi}{2\pi i} \int_{c+iT}^{c+2iT} \chi(1-s) F(s)ds + O\bigg( \frac{T^{\frac{1}{2}+\varepsilon}}{a^2} \bigg)
\end{split}\end{equation}
where $c=1+\frac{1}{\log T}$ and
$$ F(s) =  \bigg(\frac{\zeta'}{\zeta}\bigg(s+\frac{a}{\log T}\bigg)\bigg)^2\zeta(s)^2
= \sum_{n=1}^{\infty}\frac{f(n)}{n^s},$$
with
%$$f(n) =(\Lambda_a*\Lambda_a*\tau)(n), \quad \Lambda_a(m) = -\Lambda(m)m^{-a/\log T}. $$
$$f(n) =\Big[\Big(-\Lambda(\cdot)\cdot^{-a/\log T}\Big)*\Big(-\Lambda(\cdot)\cdot^{-a/\log T}\Big)*\tau \Big](n).$$
The error term $O(T^{1/2+\varepsilon}/a^2)$ comes from the horizontal contributions of the contour shift, that can be bounded exactly like in the previous proof. The factor $a^2$ in the denominator is due to the bound $\frac{\zeta'}{\zeta}(\sigma+\frac{a}{\log T}+iT)^2\ll (\log T)^2/(\sigma+\frac{a}{\log T}-\frac{1}{2})^2$ which is $\ll (\log T)^4/a^2$ for $\sigma\in(\frac{1}{2},1+\frac{1}{\log T})$.
Now we apply Lemma 2 in \cite{ConreySimpleZeros} (with $B(s)=1$), obtaining
\begin{equation}\begin{split}\label{13febr.1}
I_{|\zeta|^2}^{(2)}( a ;T) 
&=  2\pi\sum_{n\leq T/2\pi}f(n) + O\bigg(\frac{T^{\frac{1}{2}+\varepsilon}}{a^2}\bigg) .
\end{split}\end{equation}
Since
\begin{equation}\begin{split}\notag
\sum_{n\leq T/2\pi}f(n) 
%&= \sum_{n\leq T/2\pi} \sum_{d|n} \Lambda(d)d^{-a/\log T} \sum_{m|n/d}\Lambda(m)m^{-a/\log T} \tau\bigg(\frac{n}{dm}\bigg)\\
&= \sum_{d\leq T/2\pi} \Lambda(d)d^{-a/\log T} \sum_{m\leq T/2\pi d}\Lambda(m)m^{-a/\log T} \sum_{n\leq T/2\pi dm}\tau(n),
\end{split}\end{equation}
an application of Dirichlet's theorem gives
\begin{equation}\begin{split}\label{13febr.10}
\sum_{n\leq T/2\pi}f(n) = 
\bigg(1+O\bigg(\frac{1}{\log T}\bigg)\bigg) S_1 - S_2 - S_3,
\end{split}\end{equation}
with
\begin{equation}\begin{split}\notag
S_1&= \frac{T\log T}{2\pi}\sum_{d\leq T/2\pi} \frac{\Lambda(d)d^{-a/\log T}}{d} \sum_{m\leq T/2\pi d}\frac{\Lambda(m)m^{-a/\log T}}{m} \\
S_2&= \frac{T}{2\pi}\sum_{d\leq T/2\pi} \frac{\Lambda(d)d^{-a/\log T}\log d}{d} \sum_{m\leq T/2\pi d}\frac{\Lambda(m)m^{-a/\log T}}{m} \\
S_3&= \frac{T}{2\pi}\sum_{d\leq T/2\pi} \frac{\Lambda(d)d^{-a/\log T}}{d} \sum_{m\leq T/2\pi d}\frac{\Lambda(m)m^{-a/\log T}\log m}{m} .
\end{split}\end{equation}
These three sums can be estimated essentially with the same tools we used in the previous proof. 
Since (for $1\leq x\leq T$)
$$ \sum_{m\leq x}\frac{\Lambda(m)m^{-a/\log T}}{m} = \log T\frac{1-x^{-a/\log T}}{a} + O(1), $$
then
\begin{equation}\begin{split}\notag
S_1%&= \frac{T(\log T)^2}{2\pi}\sum_{d\leq T/2\pi} \frac{\Lambda(d)d^{-a/\log T}}{d} \frac{1-(\frac{T}{d})^{-a/\log T}}{a} + O(T(\log T)^2)\\
&= \frac{T(\log T)^2}{2\pi} \frac{1}{a} \sum_{d\leq T/2\pi} \frac{\Lambda(d)d^{-a/\log T}}{d} \bigg(  1-\bigg(\frac{T}{2\pi d}\bigg)^{-a/\log T}  \bigg) 
+ O\bigg( T\log T\sum_{d\leq T/2\pi} \frac{\Lambda(d)}{d}\bigg)\\
&= \frac{T(\log T)^2}{2\pi}\frac{1}{a}\sum_{d\leq T/2\pi} \frac{\Lambda(d)d^{-a/\log T}}{d} - \frac{T(\log T)^2}{2\pi}\frac{e^{-a}}{a}\sum_{d\leq T/2\pi} \frac{\Lambda(d)}{d} + O(T(\log T)^2)\\
&= \frac{T(\log T)^3}{2\pi}\frac{1}{a^2}\bigg(1-\bigg(\frac{T}{2\pi} \bigg)^{-a/\log T}  \bigg)
- \frac{T(\log T)^3}{2\pi}\frac{e^{-a}}{a} + O\bigg(\frac{T(\log T)^2}{a}\bigg) \\
&= \frac{T(\log T)^3}{2\pi}\frac{1- e^{-a}}{a^2}
- \frac{T(\log T)^3}{2\pi}\frac{e^{-a}}{a} + O\bigg(\frac{T(\log T)^2}{a}\bigg) \\
&= \frac{T(\log T)^3}{2\pi}\frac{1-e^{-a}-ae^{-a}}{a^2} + O\bigg(\frac{T(\log T)^2}{a}\bigg).
\end{split}\end{equation}
Similarly, being (for $1\leq x\leq T$)
$$ \sum_{m\leq x}\frac{\Lambda(m)m^{-a/\log T}\log m}{m} = \log T\frac{\log T - x^{-a/\log T}(a\log x+\log T)}{a^2}+O(\log T) ,$$
we have
\begin{equation}\begin{split}\notag
S_2&= \frac{T\log T}{2\pi}\sum_{d\leq T/2\pi} \frac{\Lambda(d)d^{-a/\log T}\log d}{d} \frac{1-e^{-a}d^{a/\log T}}{a} + O\bigg(T \sum_{d\leq T} \frac{\Lambda(d)\log d}{d} \bigg)\\
&=  \frac{T\log T}{2\pi}\frac{1}{a}\sum_{d\leq T/2\pi} \frac{\Lambda(d)d^{-a/\log T}\log d}{d} - \frac{T\log T}{2\pi}\frac{e^{-a}}{a}\sum_{d\leq T/2\pi} \frac{\Lambda(d)\log d}{d} + O(T(\log T)^2)\\
%&=  \frac{T(\log T)^3}{2\pi}\frac{1-ae^{-a}-e^{-a}}{a^3} - \frac{T\log T}{2\pi}\frac{e^{-a}}{a} \frac{(\log T)^2}{2} + O(T(\log T)^2)\\
&=\frac{T(\log T)^3}{2\pi}\bigg(\frac{1-ae^{-a}-e^{-a}}{a^3} - \frac{e^{-a}}{2a} \bigg) + O\bigg(\frac{T(\log T)^2}{a^2}\bigg)
\end{split}\end{equation}
and
%\begin{equation}\begin{split}\notag
%S_3%&= \frac{T\log T}{2\pi}\sum_{d\leq T/2\pi} \frac{\Lambda(d)d^{-a/\log T}}{d} \frac{\log T - (\frac{T}{d})^{-a/\log T}((a+1)\log T - \log d)}{a^2} + O(T(\log T)^2)\\
%&= \frac{T(\log T)^2}{2\pi}\frac{1}{a^2}\sum_{d\leq T/2\pi} \frac{\Lambda(d)d^{-a/\log T}}{d}  -
%\frac{T(\log T)^2}{2\pi} \frac{e^{-a}(a+1)}{a^2} \sum_{d\leq T/2\pi} \frac{\Lambda(d)}{d}  +
%\frac{T\log T}{2\pi} \frac{e^{-a}}{a^2} \sum_{d\leq T/2\pi} \frac{\Lambda(d)\log d}{d}  + O(T(\log T)^2)\\
%%&= \frac{T(\log T)^3}{2\pi}\frac{1-e^{-a}}{a^3} -
%%\frac{T(\log T)^3}{2\pi} \frac{e^{-a}(a+1)}{a^2}   +
%%\frac{T\log T}{2\pi} \frac{e^{-a}}{a^2}\frac{(\log T)^2}{2}  + O(T(\log T)^2)\\
%&= \frac{T(\log T)^3}{2\pi}\bigg(\frac{1-e^{-a}}{a^3} -
%\frac{e^{-a}(a+1)}{a^2}   +
% \frac{e^{-a}}{2a^2}\bigg) + O(T(\log T)^2).
%\end{split}\end{equation}
\begin{equation}\begin{split}\notag
S_3&= \frac{T\log T}{2\pi}\sum_{d\leq T/2\pi} \frac{\Lambda(d)d^{-a/\log T}}{d} \frac{\log T - e^{-a}d^{-a/\log T}((a+1)\log T - a\log d)}{a^2} + O\bigg(\frac{T(\log T)^2}{a}\bigg)\\
&= \frac{T(\log T)^2}{2\pi}\bigg(\frac{1}{a^2}\sum_{d\leq T/2\pi} \frac{\Lambda(d)d^{-a/\log T}}{d}  -
\frac{e^{-a}(a+1)}{a^2}\log T +
\frac{e^{-a}}{a} \frac{q\log T}{2} \bigg) + O\bigg(\frac{T(\log T)^2}{a^2}\bigg) \\
%&= \frac{T(\log T)^3}{2\pi}\frac{1-e^{-a}}{a^3} -
%\frac{T(\log T)^3}{2\pi} \frac{e^{-a}(a+1)}{a^2}   +
%\frac{T\log T}{2\pi} \frac{e^{-a}}{a^2}\frac{(\log T)^2}{2}  + O(T(\log T)^2)\\
&= \frac{T(\log T)^3}{2\pi}\bigg(\frac{1-e^{-a}}{a^3} -
\frac{e^{-a}(a+1)}{a^2}   +
 \frac{e^{-a}}{2a}\bigg) + O\bigg(\frac{T(\log T)^2}{a^2}\bigg).
\end{split}\end{equation}
Plugging these into \eqref{13febr.10}, we finally get
\begin{equation}\begin{split}\notag
\sum_{n\leq T/2\pi}f(n) 
%&=
%\frac{T(\log T)^3}{2\pi}\bigg(
%\frac{1-e^{-a}}{a^2} -\frac{e^{-a}}{a}
%-\frac{1-ae^{-a}-e^{-a}}{a^3} + \frac{e^{-a}}{2a}
%-\frac{1-e^{-a}}{a^3} + \frac{e^{-a}(a+1)}{a^2}   - \frac{e^{-a}}{2a}
%\bigg) + O(T(\log T)^2)\\
&=\frac{T(\log T)^3}{2\pi}\bigg(
\frac{a-2 +e^{-a}(a+2)}{a^3}
\bigg)+ O\bigg(\frac{T(\log T)^2}{a^2}\bigg).
\end{split}\end{equation}
In view of \eqref{13febr.1}, the proof is concluded.
\endproof

To prove Proposition \ref{lowerbound}, we then need a lower bound for the quantity 
$$\mathcal S_{|\zeta^2|}(\sigma;T) = \int_T^{2T} \sum_{\gamma,\gamma'} \frac{(\sigma-\frac{1}{2})^2}{[(\sigma-\frac{1}{2})^2+(t-\gamma)^2][(\sigma-\frac{1}{2})^2+(t-\gamma')^2]} |\zeta(\tfrac{1}{2}+it)|^2dt,$$
which we obtain in the following lemma by analyzing the \lq\lq diagonal\rq\rq \, contribution $\gamma=\gamma'$. 
We note that, even though we only need a bound from below for $\mathcal S_{|\zeta^2|}(\sigma;T)$, in the proof we proceed with asymptotic formulae as much as possible. In this way, it will be clear that the only step where the equality does not hold (see Equation \eqref{4oct.100}) is when we lower bound the double sum over zeros with its diagonal.

\begin{lemma}\label{Spesato}
Assume RH. For $\sigma=\frac{1}{2}+\frac{a}{\log T}$
with $0< a \ll 1$, 
we have
$$ \mathcal S_{|\zeta|^2}(\sigma,T) \geq  T(\log T)^3 \frac{a^2-4a+6-2e^{-a}(a+3)}{4a^3} + O\bigg(\frac{T(\log T)^{11/4}}{a}\bigg)
%+  O\bigg(\frac{(\log T)^{5}}{a^2T}\bigg) 
. $$
%In particular, for $a\to 0$ sufficiently slowly, $$ \frac{S_{|\zeta|^2}(\sigma,T)}{T(\log T)^3} \gtrsim \frac{a}{48}$$
%In particular, for $a\to 0$ sufficiently slowly, the main term above is $ \geq  \frac{a}{48} T(\log T)^3 $.
\end{lemma}

\proof
First of all, we restrict the sums over zeros in $\mathcal S_{|\zeta|^2}(\sigma,T)$ to $T\leq\gamma\leq2T$ up to a negligible error.
To show this, we write 
\begin{equation}\begin{split}\label{30may.01} 
\bigg(\sum_{\gamma}\frac{\sigma-\frac{1}{2}}{(\sigma-\frac{1}{2})^2+(t-\gamma)^2}\bigg)^2 
= &\bigg(\sum_{T\leq\gamma\leq 2T}\frac{\sigma-\frac{1}{2}}{(\sigma-\frac{1}{2})^2+(t-\gamma)^2}\bigg)^2 \\
&+ O\bigg(\sum_{\gamma}\frac{\sigma-\frac{1}{2}}{(\sigma-\frac{1}{2})^2+(t-\gamma)^2} \cdot \sum_{\gamma'\not\in[T,2T]} \frac{\sigma-\frac{1}{2}}{(\sigma-\frac{1}{2})^2+(t-\gamma')^2}\bigg)
\end{split}\end{equation}
and we integrate the above formula with respect to the measure $|\zeta(\frac{1}{2}+it)|^2dt$. To handle the error term in \eqref{30may.01} on weighted average, we appeal to the Lindel\"of Hypothesis and the following well-known bounds (see \cite[Equations (2.6) and (2.7)]{GGM}); for $T\leq t\leq  2T$: 
\begin{equation}\begin{split}\notag 
\sum_{\gamma\not\in[T,2T]}\frac{\sigma-\frac{1}{2}}{(\sigma-\frac{1}{2})^2+(t-\gamma)^2}
&\ll  \log T (\sigma-\tfrac{1}{2}) \bigg(\frac{1}{2T-t+1} + \frac{1}{t-T+1} + \frac{(\sigma-\frac{1}{2})^{-2}}{1+(t-T)^2} +  \frac{(\sigma-\frac{1}{2})^{-2}}{1+(2T-t)^2}\bigg), \\
\sum_{\gamma}\frac{\sigma-\frac{1}{2}}{(\sigma-\frac{1}{2})^2+(t-\gamma)^2} & \ll \frac{\log T}{\sigma-\frac{1}{2}}.
\end{split}\end{equation}
Doing so, we obtain 
\begin{equation}\begin{split}\label{4oct.1}\notag
\mathcal S_{|\zeta|^2}(\sigma,T) 
&=   \int_T^{2T}\bigg(\sum_{T\leq\gamma\leq 2T}\frac{\sigma-\frac{1}{2}}{(\sigma-\frac{1}{2})^2+(t-\gamma)^2} \bigg)^2 |\zeta(\tfrac{1}{2}+it)|^2dt +  O\bigg(\frac{T^{\varepsilon/2}}{(\sigma-\frac{1}{2})^2}\bigg)\\
&=\int_T^{2T}\sum_{T\leq\gamma,\gamma'\leq 2T}\frac{(\sigma-\frac{1}{2})^2|\zeta(\tfrac{1}{2}+it)|^2}{[(\sigma-\frac{1}{2})^2+(t-\gamma)^2][(\sigma-\frac{1}{2})^2+(t-\gamma')^2]}  dt  + O\bigg( \frac{T^{\varepsilon}}{a^2}\bigg).
\end{split}\end{equation}
for any fixed $\varepsilon>0$. We consider the diagonal $\gamma=\gamma'$ in the main term above, denoting
\begin{equation}\label{4oct.2}
\mathcal S_{|\zeta|^2}^{\mathcal D}(\sigma,T) 
:=   \int_T^{2T}\sum_{T\leq\gamma\leq 2T}\frac{(\sigma-\frac{1}{2})^2}{[(\sigma-\frac{1}{2})^2+(t-\gamma)^2]^2}  |\zeta(\tfrac{1}{2}+it)|^2dt .
\end{equation}
Trivially
%\begin{equation}\label{4oct.100}
%\mathcal S_{|\zeta|^2}(\sigma,T) 
%\geq \mathcal S_{|\zeta|^2}^{\mathcal D}(\sigma,T) +O \bigg(\frac{T^{\varepsilon}}{a^2} \bigg).
%\end{equation}
\begin{equation}\label{4oct.100}
\mathcal S_{|\zeta|^2}(\sigma,T) +O \bigg(\frac{T^{\varepsilon}}{a^2} \bigg)
= \int_T^{2T}\sum_{T\leq\gamma,\gamma'\leq 2T}\frac{(\sigma-\frac{1}{2})^2|\zeta(\tfrac{1}{2}+it)|^2}{[(\sigma-\frac{1}{2})^2+(t-\gamma)^2][(\sigma-\frac{1}{2})^2+(t-\gamma')^2]}  dt 
\geq \mathcal S_{|\zeta|^2}^{\mathcal D}(\sigma,T) .
\end{equation}
Therefore, it suffices to lower bound $ \mathcal S_{|\zeta|^2}^{\mathcal D}(\sigma,T)$. The contribution from $|t-\gamma|\geq 1$ in \eqref{4oct.2} is $\ll T$, since
\begin{equation}\begin{split}\notag
\sum_{T\leq\gamma\leq 2T} \int_{\substack{t\in[T,2T]:\\|t-\gamma|\geq 1}}\frac{(\sigma-\frac{1}{2})^2}{[(\sigma-\frac{1}{2})^2+(t-\gamma)^2]^2}  |\zeta(\tfrac{1}{2}+it)|^2dt &
\ll \int_T^{2T} \sum_{\substack{T\leq\gamma\leq 2T \\ |\gamma-t|\geq 1}} \frac{(\sigma-\frac{1}{2})^2}{(t-\gamma)^4}  |\zeta(\tfrac{1}{2}+it)|^2dt \\
&\ll \log T(\sigma-\tfrac{1}{2})^2 \int_{T}^{2T}  |\zeta(\tfrac{1}{2}+it)|^2dt 
%\ll  \frac{a^2}{\log T} T\log T  \ll a^2T.
\ll T.
\end{split}\end{equation}
Hence, with the change of variable $t-\gamma = x$, we have
\begin{equation}\begin{split}\notag
\mathcal S_{|\zeta|^2}^{\mathcal D}(\sigma,T) 
&= \sum_{T\leq\gamma\leq 2T}  \int_{\substack{t\in[T,2T]:\\|t-\gamma|\leq 1}}\frac{(\sigma-\frac{1}{2})^2}{[(\sigma-\frac{1}{2})^2+(t-\gamma)^2]^2}  |\zeta(\tfrac{1}{2}+it)|^2dt + O(T)\\
&= \int_{-1}^{1}\frac{(\sigma-\frac{1}{2})^2}{[(\sigma-\frac{1}{2})^2+x^2]^2}  \sum_{T\leq\gamma\leq 2T}\left|\zeta\left(\tfrac{1}{2}+i\gamma+ix\right)\right|^2dx + O(T).
%&= \int_{-1}^{1}\frac{(\sigma-\frac{1}{2})^2}{[(\sigma-\frac{1}{2})^2+x^2]^2}  \sum_{T\leq\gamma\leq 2T}\bigg|\zeta\bigg(\frac{1}{2}+i\gamma+i\frac{\frac{x}{2\pi}\log T}{\frac{1}{2\pi}\log T}\bigg)\bigg|^2dx + O(T).
\end{split}\end{equation}
The sum over zeros can be evaluated by Gonek's asymptotic formula for the second discrete moment of zeta \cite[Corollary 2]{Gonek} (see \cite[Theorem 2]{GonekBook} for a little bit more of uniformity), getting
\begin{equation}\begin{split}\notag
\mathcal S_{|\zeta|^2}^{\mathcal D}(\sigma,T) 
= \frac{T(\log T)^2}{2\pi}\int_{-1}^{1}\frac{(\sigma-\frac{1}{2})^2}{[(\sigma-\frac{1}{2})^2+x^2]^2}  
\bigg(1-&\frac{\sin^2(\frac{x}{2}\log T)}{(\frac{x}{2}\log T)^2}\bigg) 
dx \\ 
&+ O\bigg(T(\log T)^{7/4}\int_{-1}^{1}\frac{(\sigma-\frac{1}{2})^2}{[(\sigma-\frac{1}{2})^2+x^2]^2}\bigg) + O(T).
\end{split}\end{equation}
In view of \cite[p. 114, fourth display]{GGM}, the remaining integral in the error term can be bounded by $\ll (\sigma-\frac{1}{2})^{-1}$. Therefore, with the change of variable $y=x\log T$ in the main term above, we get
\begin{equation}\begin{split}\notag
\mathcal S_{|\zeta|^2}^{\mathcal D}(\sigma,T) 
&= \frac{T(\log T)^3}{2\pi}\int_{-\log T}^{\log T}\frac{a^2}{(a^2+y^2)^2}  
\bigg(1-\frac{\sin^2(\frac{y}{2})}{(\frac{y}{2})^2}\bigg) 
dy + O\bigg(\frac{T(\log T)^{11/4}}{a}\bigg) .
\end{split}\end{equation}
Since
\begin{equation}\notag
\int_{|y|>\log T}\frac{a^2}{(a^2+y^2)^2} \bigg(1-\frac{\sin^2(\frac{y}{2})}{(\frac{y}{2})^2}\bigg) dy
\ll  \int_{|y|>\log T}\frac{1}{y^4}dy \ll \frac{1}{(\log T)^3}
\end{equation} we have
\begin{equation}\label{4oct.3}
\mathcal S_{|\zeta|^2}^{\mathcal D}(\sigma,T) 
= \frac{T(\log T)^3}{2\pi}\int_{-\infty}^{+\infty}\frac{a^2}{(a^2+y^2)^2}  
\bigg(1-\frac{\sin^2(\frac{y}{2})}{(\frac{y}{2})^2}\bigg) 
dy + O\bigg(\frac{T(\log T)^{11/4}}{a}\bigg) .
\end{equation}
We note that, for $a>0$, by the residue theorem we get
\begin{equation}\begin{split}\label{intergraleresiduipesato}\notag
\frac{1}{2\pi}\int_{-\infty}^{+\infty}\frac{a^2}{(a^2+y^2)^2}  \bigg(1-\frac{\sin^2(\frac{y}{2})}{(\frac{y}{2})^2}\bigg)dy 
&= \frac{a^2+6\sinh(a)-2a(\cosh(a)+2)+2a\sinh(a)-6\cosh(a)+6}{4a^3}\\
%&=2\pi \frac{a^2-2a(e^{-a}+2)-6(e^{-a}-1)}{4a^3}
&=\frac{a^2-4a+6-2e^{-a}(a+3)}{4a^3}.
\end{split}\end{equation}
Plugging the above in \eqref{4oct.3}, we have
$$ \mathcal S_{|\zeta|^2}^{\mathcal D}(\sigma,T) 
= T(\log T)^3\frac{a^2-4a+6-2e^{-a}(a+3)}{4a^3}  + O\bigg(\frac{T(\log T)^{11/4}}{a}\bigg) $$
that ends the proof in view of \eqref{4oct.100}.
%Hence, by \eqref{4oct.100} 
%\begin{equation}\begin{split}\notag
%S_{|\zeta|^2}(\sigma,T) 
%&\geq S_{|\zeta|^2}^{\mathcal D}(\sigma,T) + O((\log T)^3)+ \bigg(\frac{ (\log T)^5}{a^2T} \bigg)\\
%&=T(\log T)^3 \frac{a^2-4a+6-2e^{-a}(a+3)}{4a^3}  +  O\bigg(\frac{T(\log T)^{11/4}}{a}\bigg) +  O\bigg(\frac{(\log T)^{5}}{a^2T}\bigg).
%\end{split}\end{equation} 
\endproof

%Proposition \ref{lowerbound} trivially follows from Lemma \ref{lemmaSpezzamentoI} and \ref{Spesato}  above, as (for $\frac{\log T}{T}\ll a \ll 1$)
%\begin{equation}\begin{split}\notag 
%I_{|\zeta|^2}( a ;T) &
%= T(\log T)^3 \frac{(a-2)(2e^{-a}-2-a^2+2a)}{2 a^3} + 2 \mathcal S_{|\zeta|^2}(\sigma,T) + O(T(\log T)^2)\\
%&\geq T(\log T)^3\bigg( \frac{(a-2)(2e^{-a}-2-a^2+2a)}{2 a^3} +  \frac{a^2-4a+6-2e^{-a}(a+3)}{2a^3}\bigg) \\
%&\hspace{5cm}+ O(a^2T) + O\bigg(\frac{T(\log T)^{11/4}}{a}\bigg) + O(T(\log T)^2) \\
%&= T(\log T)^3 \frac{ - a^3 +5 a^2 -10 a - 10 (e^{-a}-1)}{2 a^3} \\
%&\hspace{5cm}+ O(a^2T) + O\bigg(\frac{T(\log T)^{11/4}}{a}\bigg) + O(T(\log T)^2) .
%\end{split}\end{equation}
%\begin{equation}\begin{split}\notag 
%I_{|\zeta|^2}( a ;T) &
%= T(\log T)^3 \frac{(a-2)(2e^{-a}-2-a^2+2a)}{2 a^3} + 2 \mathcal S_{|\zeta|^2}\bigg(\frac{1}{2}+\frac{a}{\log T},T\bigg) + O(T(\log T)^2)\\
%&\geq T(\log T)^3\bigg( \frac{(a-2)(2e^{-a}-2-a^2+2a)}{2 a^3} +  \frac{a^2-4a+6-2e^{-a}(a+3)}{2a^3}\bigg) + O\bigg(\frac{T(\log T)^{11/4}}{a}\bigg)  \\
%&= T(\log T)^3 \frac{ - a^3 +5 a^2 -10 a - 10 (e^{-a}-1)}{2 a^3} + O\bigg(\frac{T(\log T)^{11/4}}{a}\bigg)  .
%\end{split}\end{equation}
%In the range  $(\log T)^{-1/4+\varepsilon}\ll a \ll 1$ the error terms are $\ll (\log T)^{3-\varepsilon}$; finally for $a=o(1)$
%$$ \frac{ - a^3 +5 a^2 -10 a - 10 (e^{-a}-1)}{2 a^3} 
%= \frac{1}{3} - \frac{5 a}{24} + \frac{ a^2}{24} + O(a^3) 
%\geq \frac{1}{3} - \frac{5 a}{24}.$$

\section{Upper bound}\label{SectionUpperBound}

First of all we notice that, by \cite[Equation (2.6)]{GGM}, one can easily see that under RH
\begin{equation}\notag 
\mathcal S_{|\zeta|^2}(\sigma,T) 
\ll \int_{T}^{2T} \bigg(\frac{\log t}{\sigma-\frac{1}{2}}\bigg)^2 |\zeta(\tfrac{1}{2}+it)|^2 dt 
\ll \frac{(\log T)^4}{a^2}\int_{T}^{2T}|\zeta(\tfrac{1}{2}+it)|^2dt 
\ll \frac{T(\log T)^5}{a^2},
\end{equation}
which immediately implies (together with Lemma \ref{lemmaSpezzamentoI})
$$ I_{|\zeta|^2}(a;T) \ll  \frac{T(\log T)^5}{a^2}. $$

Assuming an hypothesis on gaps between zeros of zeta (see \cite{FanGe}, Equation (ES $2K$)), the power of $\log T$ can be reduced to 4, by applying the Cauchy-Schwarz inequality. Namely, applying e.g. \cite[Theorem 1.1]{FanGe} with $K=1$, for $a=o(1)$ as $T\to\infty$, we have
\begin{equation}\begin{split}\notag 
I_{|\zeta|^2}(a,T) 
&\ll \sqrt{\int_{T}^{2T} \bigg|\frac{\zeta'}{\zeta}\bigg(\frac{1}{2}+\frac{a}{\log T}+it\bigg)\bigg|^4 dt} \times \sqrt{\int_{T}^{2T} |\zeta(\tfrac{1}{2}+it)|^4 dt}\\
&\ll \sqrt{\frac{T(\log T)^4}{a^3}} \times \sqrt{T(\log T)^4}
\ll \frac{T(\log T)^4}{a^{3/2}}.
\end{split}\end{equation}

The main result of this section is an unconditional (i.e. no assumptions on gaps, RH is always assumed here) upper bound for $I_{|\zeta|^2}(a;T)$. 
%The power of $\log T$ in our bound is optimal, unlike the dependence on $a$.
The power of $\log T$ in our bound is optimal. Moreover, due to our weight $|\zeta(\frac{1}{2}+it)|^2$ that cancels the denominator of 
%$|\frac{\zeta'}{\zeta}(\frac{1}{2}+\frac{a}{\log T}+it)|^2$
$|\zeta'/\zeta(\frac{1}{2}+\frac{a}{\log T}+it)|^2$
as  $a$ approaches 0, we have no loss in $a$ when $a=o(1)$. However, for a fixed $a$, the constant in our upper bound is not sharp.

\begin{prop}\label{upperbound}
Assume RH. For $0< a \ll 1$, as $T\to\infty$, we have
$$I_{|\zeta|^2}(a;T) \leq T(\log T)^3\frac{e^a-e^{-a}(2a^2+2a+1)}{4a^3} + O(T(\log T)^2).$$
Moreover, if $a=o(1)$, then
$$I_{|\zeta|^2}(a;T) \leq T(\log T)^3\bigg(\frac{1}{3} - \frac{a}{6} + \frac{a^2}{15}\bigg) + O(T(\log T)^2).$$
\end{prop}

\proof
Assuming RH, Soundararajan \cite{Sound} (Equation (5) and above) proved that, for $\sigma>\frac{1}{2}$ and $t\in[T,2T]$ with $t\neq\gamma$, we have
\begin{equation}\begin{split}\notag 
\log|\zeta(\tfrac{1}{2}+it)|& - \log|\zeta(\sigma+it)|  
= (\sigma-\tfrac{1}{2})\bigg(\frac{\log T}{2}+O(1)\bigg) -\frac{1}{2}\sum_{\gamma}\log\frac{(\sigma-\frac{1}{2})^2+(t-\gamma)^2}{(t-\gamma)^2}\\
&\leq (\sigma-\tfrac{1}{2})\bigg(\frac{\log T}{2}+O(1)-\frac{1}{2}\sum_\gamma\frac{\sigma-\frac{1}{2}}{(\sigma-\frac{1}{2})^2+(t-\gamma)^2}\bigg)
\leq (\sigma-\tfrac{1}{2})\bigg(\frac{\log T}{2}+O(1)\bigg).
\end{split}\end{equation}
Therefore
%\begin{equation}\begin{split}\notag 
% \log\bigg|\zeta\bigg(\frac{1}{2}+it\bigg)\bigg| - \log\bigg|\zeta\bigg(\frac{1}{2}+\frac{a}{\log T}+it\bigg)\bigg|
%&\leq \frac{a}{2}+O\bigg(\frac{a}{\log T}\bigg)-\frac{1}{2}\sum_\gamma\frac{(\frac{a}{\log T})^2}{(\frac{a}{\log T})^2+(t-\gamma)^2}\\
%&\leq \frac{a}{2}+O\bigg(\frac{a}{\log T}\bigg) 
%\end{split}\end{equation}
%and then
\begin{equation}\begin{split}\notag
\bigg|\frac{\zeta(\frac{1}{2}+it)}{\zeta(\frac{1}{2}+\frac{a}{\log T}+it)}\bigg|
&= \exp\bigg( \log\bigg|\zeta\bigg(\frac{1}{2}+it\bigg)\bigg| - \log\bigg|\zeta\bigg(\frac{1}{2}+\frac{a}{\log T}+it\bigg)\bigg| \bigg)
%&\leq \exp\bigg( \frac{a}{2}+O\bigg(\frac{a}{\log T}\bigg)-\frac{1}{2}\sum_\gamma\frac{(\frac{a}{\log T})^2}{(\frac{a}{\log T})^2+(t-\gamma)^2} \bigg)
\leq \exp\bigg( \frac{a}{2}+O\bigg(\frac{a}{\log T}\bigg)\bigg).
\end{split}\end{equation}
From the above we deduce
\begin{equation}\begin{split}\notag 
I_{|\zeta|^2}(a,T) 
&= \int_{T}^{2T} \bigg|\frac{\zeta(\frac{1}{2}+it)}{\zeta(\frac{1}{2}+\frac{a}{\log T}+it)}\bigg|^2\bigg|\zeta'\bigg(\frac{1}{2}+\frac{a}{\log T}+it\bigg)\bigg|^2 dt\\
%&\leq \int_{T}^{2T} \exp\bigg( a+O\bigg(\frac{a}{\log T}\bigg)\bigg)\bigg|\zeta'\bigg(\frac{1}{2}+\frac{a}{\log T}+it\bigg)\bigg|^2 dt \\
&\leq e^{ a}\bigg(1+O\bigg(\frac{a}{\log T}\bigg)\bigg)\int_{T}^{2T} \bigg|\zeta'\bigg(\frac{1}{2}+\frac{a}{\log T}+it\bigg)\bigg|^2 dt.
\end{split}\end{equation}
Ingham's work \cite{Ingham} yields
\begin{equation}\begin{split}\notag 
\int_{T}^{2T} \bigg|\zeta'\bigg(\frac{1}{2}+\frac{a}{\log T}+it\bigg)\bigg|^2 dt = T(\log T)^3\frac{2-e^{-2a}(4a^2+4a+2)}{(2a)^3} + O(T(\log T)^2),
\end{split}\end{equation}
so 
\begin{equation}\begin{split}\notag 
\frac{I_{|\zeta|^2}(a,T)}{T(\log T)^3} 
%&\leq e^{ a}\bigg(1+O\bigg(\frac{a}{\log T}\bigg)\bigg)\frac{e^{-2a}(-4a^2-4a-2)+2}{(2a)^3} + O\bigg(\frac{1}{\log T}\bigg)\\
&\leq \frac{e^a-e^{-a}(2a^2+2a+1)}{4a^3} + O\bigg(\frac{1}{\log T}\bigg).
\end{split}\end{equation}
Moreover, if $a=o(1)$, then
$$ \frac{e^a-e^{-a}(2a^2+2a+1)}{4a^3} 
= \frac{1}{3} - \frac{a}{6} + \frac{a^2}{15}-\frac{a^3}{60} + O(a^4)
\leq  \frac{1}{3} - \frac{a}{6} + \frac{a^2}{15}.$$
\endproof

{\small

\end{document}